\documentclass[12pt]{article}
\usepackage{amsmath,amsfonts,amscd,a4wide,theorem}
\theoremheaderfont{\bfseries\rmfamily\upshape}
                \newtheorem{theorem}{Theorem}

                \newtheorem{lemma}{Lemma}

{\theorembodyfont{\rmfamily}

}

\newcommand{\qed}{\text{\rule{.4em}{1.7ex}\hspace{.6em}}}
\newenvironment{proof}[1][]{\noindent {\bf Proof#1:\ }}
        {\hspace*{.1em}\hfill\qed\bigskip \noindent}
\newcounter{rom}
\renewcommand{\therom}{(\roman{rom})}
{\end{list}}
\title{Lagrangian submanifolds attaining equality in the improved Chen's inequality}
\begin{document}
\author{J. Bolton and  L. Vrancken}


\maketitle

\sloppy

\begin{abstract}\noindent
In \cite{O} Oprea gave an improved version of Chen's inequality for Lagrangian submanifolds of $\mathbb CP^n(4)$. For minimal
submanifolds this inequality coincides with a previous version proved in \cite{CDVV1}. We consider here those non
minimal $3$-dimensional Lagrangian submanifolds in $\mathbb CP^3 (4)$ attaining at all points equality in the improved Chen inequality.
 We show how all such submanifolds may be obtained starting from a minimal Lagrangian surface in $\mathbb CP^2(4)$.
\end{abstract}

{{\bfseries Key words}:  {\em Lagrangian submanifold, complex projective space, Chen inequality}.

{\bfseries Subject class: } 53B25, 53B20.}

\footnotetext{This work was done during a research visit of the second author at Durham University, both authors are grateful for the support of the universities of Durham and Valenciennes}

\section{Introduction}
In the early nineties Chen \cite{C} introduced a new invariant, called $\delta_M$, for a Riemannian manifold $M$. Specifically,  $\delta_M:M\to {\mathbb R}$ is given by:
 $$\delta_M (p) =\tau(p) - (\inf K)(p),$$
where
$(\inf K)(p) = \inf \big\{ K(\pi)\, | \,\pi \text{ is a 2-dimensional
subspace of } T_p M \big\}$,
with $K(\pi)$ being the sectional curvature of $\pi$, and
$\tau(p)=\sum_{i<j}K(e_i\wedge e_j)$ denotes the scalar curvature
defined in terms of an orthonormal basis $\{e_1,\ldots,e_n\}$ of the
tangent space $T_p M$ of $M$ at $p$.  In the same paper, he discovered, for submanifolds of real space forms, an inequality relating this invariant with the length of the mean curvature vector $H$.  A similar inequality was proved in \cite{CDVV1} and \cite{CDVV2} for n-dimensional Lagrangian submanifolds of a complex space form $\tilde M^n(4c)$ of constant holomorphic sectional curvature $4 c$. Indeed, it was shown that
\begin{equation} \label{oldcheneq}
\delta_M \le \tfrac{(n-2)(n+1)}{2} \tfrac{c}{4}  + \tfrac{n^2}{2} \tfrac{n-2}{n-1} \Vert H \Vert^2 \, .
\end{equation}
Note that, for $n=2$, both sides of the above inequality are zero.

Let  $\mathbb CP^n(4)$ denote complex projective n-space of constant holomorphic sectional curvature 4. For $n\ge 3$, Lagrangian submanifolds of $\mathbb CP^n(4)$ attaining at every point  equality in \eqref{oldcheneq}  were studied in, amongst others, \cite{CDVV1}, \cite{CDVV2}, \cite{BSVW} and \cite{BSV}. In particular, in \cite{CDVV1} and  \cite{CDVV2}, it was shown that such submanifolds are minimal, and a complete classification was obtained of $3$-dimensional Lagrangian submanifolds of  $\mathbb CP^n(4)$  attaining at each point   equality in \eqref{oldcheneq}. Such submanifolds are obtained starting from minimal surfaces with ellipse of curvature a circle in the unit $5$-sphere.

However, Oprea \cite{O} has recently shown that the inequality  \eqref{oldcheneq} is not optimal, and, for $n\ge 3$ can be improved to
\begin{equation} \label{newcheneq}
\delta_M \le \tfrac{(n-2)(n+1)}{2} \tfrac{c}{4}  + \tfrac{n^2}{2} \tfrac{2 n-3}{2 n+3} \Vert H \Vert^2 \, .
\end{equation}

This explains why a Lagrangian submanifold of  $\mathbb CP^n(4)$ attaining  at every point  equality in \eqref{oldcheneq} must be minimal, since both inequalities coincide in this case.

\section{Classification}

Let $M$ be a Lagrangian submanifold of  $\mathbb CP^n(4)$. A careful analysis of Oprea's arguments shows that equality in \eqref{newcheneq} is obtained at a point $p\in M$ if and only if there exists an orthonormal basis $\{e_1,e_2,\dots,e_n\}$ of the tangent space $T_pM$ such that the symmetric cubic form $C$ \cite{} on $M$ constructed using the second fundamental form $h$,  the complex structure $J$ and the Riemannian metric $< \ , \ >$ on $\mathbb CP^n(4)$ given by
$$
C(X,Y,Z)= <h(X,Y),JZ>,
$$
has the following form,
\begin{align}%
<h(e_2,e_2),Je_2> &=  -<h(e_3,e_3),Je_2>\label{C1}\\
4<h(e_2,e_2),Je_1>=4 <h(e_3,e_3),Je_1>&= <h(e_1,e_1),Je_1>=3
<h(e_j,e_j),Je_1>,\label{C2}
\end{align}
where $j \in \{4,\dots, n\}$, and all other components of C are zero unless they can be obtained from the above using the symmetric nature of $C$.

In this paper we show that the inequality \eqref{newcheneq} is optimal, and we show how to construct all 3-dimensional non-minimal Lagrangian submanifolds of $\mathbb CP^3(4)$ which attain everywhere equality (the classification in the minimal case having been found in  \cite{CDVV1} and  \cite{CDVV2}).

We now assume that $M$ is a non-minimal Lagrangian submanifold of  $\mathbb CP^3(4)$ attaining at all points equality in the improved Chen inequality \eqref{newcheneq}. Then $C$ satisfies \eqref{C1} and \eqref{C2} at all points. Thus, using the notation and terminology of  \cite{RV},  $M$ is a non-minimal submanifold of Type 2 with the additional condition that $\lambda_1 =4 \lambda_2 \ne 0$, where  $\lambda_1= <h(e_1,e_1),Je_1>$ and  $\lambda_2= <h(e_2,e_2),Je_1>$. We have chosen the above orthonormal basis  $\{e_1,e_2,e_3\}$ so that the notation agres with  \cite{RV}, and, in particular, the plane for which the minimal sectional curvature is attained is that spanned by $e_2$ and $e_3$.

Since $M$ is Lagrangian, there is a horizontal lift $E_0:M\to S^7(1)\subset \mathbb R^8=\mathbb C^4$ to the unit 7-sphere \cite{R}, and if $dE_0$ denotes the derivative of $E_0$, we put $E_j=dE_0(e_j)$, for $j=1,2,3$. We will often identify a point of $M$ with its image under $E_0$.

It follows from  \cite{RV} that, for some suitable function $b_1$,
\begin{align}
&D_{E_1} E_1 = 4 \lambda_2 i E_1-E_0,\label{de1}\\
&D_{E_j} E_1 =(b_1+i \lambda_2) E_j,\qquad j=2,3,\label{dej}
\end{align}
where $D$ denotes the standard flat covariant derivative on $\mathbb
C^4$. We also get from (41), (42), (50) and (51) of \cite{RV} that
the functions $\lambda_2$ and $b_1$ have zero derivative with
respect to $E_2$ and $E_3$ and from (40) and (46) of \cite{RV} that
their derivatives in the $E_1$ direction are given by
\begin{align}
&E_1(\lambda_2) = 2 \lambda_2 b_1,\label{dl2}\\
&E_1(b_1) = -(1 +b_1^2 +3 \lambda_2^2).\label{db1}
\end{align}

The following lemma is immediate from \eqref{de1} and  \eqref{dej}.

\begin{lemma} \label{lb}The brackets $[E_1,E_2]$,  $[E_1,E_3]$, $[E_2,E_3]$ are linear combinations of $E_2$ and $E_3$.
\end{lemma}

In \cite{RV}, submanifolds of the type we are considering are
divided into $3$ further subcases depending on the relative values
of $a= <h(e_2,e_2),Je_2>$ and $\lambda_2$. One of these cases is
easy to deal with, namely that in which $a \ne 0$ but $a^2-2
\lambda_2^2 =0$ In this case, it  follows from equations (33) - (45)
of \cite{RV} that $b_1=0$ which contradicts \eqref{db1}. Hence this
case cannot occur.

We now consider the other two cases, namely those where  $a=0$, or
both $a$ and $a^2-2 \lambda_2^2$ are non-zero. We introduce a
function $\theta$ defined locally on $M$ having zero derivative with
respect to $E_2$ and $E_3$ and satisfying $E_1(\theta)=-\lambda_2$.
It follows from Lemma \ref{lb} that the integrability conditions of
the this system for $\theta$ are satisfied, and hence such a
function $\theta$ exists.

We now consider the maps into $S^7(1)$ given by
\begin{align}
&V=(-(b_1+ i \lambda_2) E_0 +E_1)/\sqrt{1+b_1^2 +\lambda_2^2},\label{V}\\
&W=e^{i\theta}(E_0 -(-b_1+i\lambda_2)E_1)/\sqrt{1+b_1^2 +\lambda_2^2}.\label{W}
\end{align}
It follows easily that $D_{E_2}V=D_{E_3}V =0$ and $D_{E_1} V = 3 \lambda_2 i V$. This implies that $V$ is contained in the unit circle of a complex plane $\mathbb C$, and, taking $t$ as the standard parameter along this circle, we also have that $E_1(t) = 3 \lambda_2$. Hence, after applying a translation if necessary, we may assume that $\theta= -t/3$.

\begin{lemma} \label{LW} The map $W$ describes a minimal horizontal surface in the unit sphere $S^5(1)$ of the orthogonal complement in $\mathbb C^4$ of the complex plane containing $V$.
\end{lemma}

\begin{proof}It is clear that $W$ is othogonal to $V$ and $iV$, so the image of $W$ is contained in the
 indicated $S^5(1)$. We now use arguments similar to those employed for $V$ above to complete the proof.
 In fact,
\begin{align*}
&D_{E_1} W=0,\\
&D_{E_j} W=\sqrt{1+b_1^2 +\lambda_2^2} e^{i\theta} E_j, \quad j=2,3,\\
&D_{E_2}(D_{E_2}W)=b_3 D_{E_3} W+ i a D_{E_2} W -(1+b_1^2+\lambda_2^2) W, \\
&D_{E_2}(D_{E_3}W)=-b_3 D_{E_2} W -i a D_{E_3} W, \\
&D_{E_3}(D_{E_2}W)=c_2 D_{E_3} W  -i a D_{E_3} W,\\
&D_{E_3}(D_{E_3}W)=-c_2 D_{E_2} W-i a D_{E_2} W   -(1+b_1^2+\lambda_2^2) W, \\
\end{align*}
from which the proof of the lemma quickly follows.
\end{proof}

We can now state  and prove our classification theorem.

\begin{theorem} Let M be a 3 dimensional non-minimal Lagrangian submanifold of  $\mathbb CP^3(4)$ which attains equality
at every point in Oprea's improvement \eqref{newcheneq} of Chen inequality. Then there is
a minimal Lagrangian surface $\tilde W$ in $\mathbb CP^2(4)$ such that $M$ can be locally written as $[E_0]$ where
\begin{equation*}
E_0 =  \frac{e^{it/3}}{\sqrt{1+b_1^2 +\lambda_2^2}}(0,W)+
\frac{(-b_1 +i\lambda_2)}{\sqrt{1+b_1^2
+\lambda_2^2}}(e^{it},0,0,0),
\end{equation*}
where $b_1$ and $\lambda_2$ are solutions of the following  system of ordinary differential equations:
\begin{equation}\label{db1dl2}
\frac{db_1}{dt}= -\frac{1+3\lambda_2^2+b_1^2}{3 \lambda_2}, \qquad
\frac{d\lambda_2}{dt}= \frac{2}{3} b_1.
\end{equation}
Conversely any $3$ dimensional Lagrangian submanifold obtained in
this way attains equality at each point in \eqref{newcheneq}.
\end{theorem}

\begin{proof}
By \cite{R}, minimal horizontal surfaces in $S^5(1)$ correspond to
minimal Lagrangian surfaces in $\mathbb CP^2(4)$. Solving \eqref{V}
and \eqref{W} for $E_0$, we find that, after applying a suitable
element of $SU(4)$, the original immersion is the projection onto
$\mathbb CP^3(4)$ of the map $E_0$ given above, where, from
\eqref{dl2} and   \eqref{db1}, $b_1$ and $\lambda_2$ are solutions
of the system \eqref{db1dl2}.
Conversely, it is clear that any submanifold obtained in this way
satisfies \eqref{C1} and \eqref{C2} and hence attains equality in
\eqref{newcheneq} at each point.
\end{proof}

\noindent{\bf Remarks} (i) It is clear that  $\lambda_2(1+\lambda_2^2+b_1^2)$ is a first integrand of the system \eqref{db1dl2}.

\noindent (ii) An alternative method of proof would be to apply immediately Theorem 7 or Theorem 9 of \cite{RV}. However the result in that case would have been less explicit.

\noindent (iii) Lagrangian immersions into ${\mathbb C}P^n(4)$ constructed from a curve in $S^3(1)$ and a lower dimensional Lagrangian immersion have been studied in \cite{CLU}.

\noindent {\em J. Bolton, Dept of Mathematical Sciences, University of Durham,
Durham DH1 3LE, UK. \quad  E-mail: john.bolton@dur.ac.uk

\noindent L. Vrancken,   LAMATH, ISTV2, Universit\'e de Valenciennes,
Campus du Mont Houy, 59313 Valenciennes Cedex 9, France. \quad E-mail:
luc.vrancken@univ-valenciennes.fr}

\end{document}